\documentclass[11pt]{amsart} 

\usepackage[utf8]{inputenc} 
\usepackage{amsfonts}
\usepackage{amsopn}
\usepackage{hyperref}
\usepackage{amsthm}
\usepackage{amsmath}
\usepackage{tikz-cd}
\usepackage{amssymb}
\usepackage{anyfontsize}

\newtheorem{theorem}{Theorem}
\newtheorem{problem}[theorem]{Problem}
\newtheorem{proposition}[theorem]{Proposition}

\newtheorem{lemma}[theorem]{Lemma}

\newtheorem{remark}[theorem]{Remark}
\newtheorem{definition}[theorem]{Definition}
\newtheorem{example}[theorem]{Example}
\numberwithin{theorem}{section}
\newcommand\pa{{\mathfrak a}}
\newcommand\pb{{\mathfrak b}}

\newcommand{\bu}{\pmb{u}}
\newcommand{\bv}{\pmb{v}}
\newcommand{\bw}{\pmb{w}}
\newcommand{\bx}{\pmb{x}}
\newcommand{\by}{\pmb{y}}
\newcommand{\bze}{\pmb{0}}

\newcommand{\bal}{\pmb{\alpha}}
\newcommand{\bbe}{\pmb{\beta}}
\newcommand{\bde}{\pmb{\delta}}
\newcommand{\bga}{\pmb{\gamma}}
\newcommand{\cF}{\mathcal{F}}

\usepackage{geometry} 
\usepackage{graphicx} 

\title{Local resolution of the second log-Fitting ideal}
\author{Yueting Jiang}
\address{Universit\'e Paris Cit\'e and Sorbonne Universit\'e, UFR de Math\'ematiques, Institut de Math\'ematiques de Jussieu-Paris Rive Gauche, UMR7586, F-75013 Paris, France.}
\email{yueting.jiang@imj-prg.fr}

\begin{document}

\begin{abstract}
    Logarithmic Fitting ideals encode the complexity of the corresponding morphisms in a natural and intrinsic way, and their principalisation offers a new perspective on problems such as the resolution of cotangent sheaves and the monomialisation of morphisms. In this paper, we give a local resolution theorem for log-rank adapted morphisms: after finite sequences of smooth local blow-ups, the first two log-Fitting ideals $\cF_{n-1}(\sigma)$ and $\cF_{n-2}(\sigma)$ become principal monomial ideals.
\end{abstract}

\maketitle

\section{Introduction}
Resolution and monomialisation problems for morphisms occupy a central position in modern birational geometry and singularity theory. Since the pioneering works of Zariski \cite{Zar39} and Hironaka \cite{hironaka1964resolution}, a large body of literature has been devoted to understanding how singular geometric or differential structures simplify after suitable birational modifications; see, for instance, \cite{cutkosky1999local,cutkosky2002monomialization,cutkosky2004monomialization,abramovich2002torification}. Beyond the classical problem of resolving singularities of varieties, one is naturally led to study the behaviour of intrinsic differential objects attached to morphisms, such as cotangent sheaves, logarithmic differential forms and Jacobian ideals \cite{da2017resolution,bbmono}.

One particularly fruitful viewpoint consists in encoding the complexity of a morphism through its logarithmic Fitting ideals. These ideals arise naturally from the logarithmic Jacobian matrix of a morphism of pairs and provide an intrinsic measure of the failure of the morphism to admit a monomial form. The use of Fitting ideals in the study of
morphisms goes back at least to the work
of Teissier \cite{Teissier1977TheHO}, and more recently they have played an important role in the study of cotangent sheaves and monomialisation problems; see for instance \cite{da2017resolution,bbmono}. In particular, the authors of \cite{da2017resolution} showed that, for birational morphisms in dimension three, the principal monomiality of logarithmic Fitting ideals is closely related to the resolution of cotangent sheaves.

Beyond its intrinsic interest, the resolution of cotangent sheaves is also motivated by the study of \(L^2\)-cohomology on the smooth locus of singular complex varieties. Building on Cheeger's foundational work on the Hodge theory of singular spaces \cite{cheeger1979spectral,cheeger1980hodge}, Cheeger, Goresky and MacPherson conjectured a relationship between \(L^2\)-cohomology and middle-perversity intersection cohomology \cite{cheeger1982l2}. The work of Hsiang and Pati on normal algebraic surfaces \cite{HP85} illustrates the role played by resolutions admitting particularly simple local descriptions of differential forms; see \cite{da2017resolution} and the references therein.

The present article is motivated by the problem of understanding the extent to which these ideas remain valid in arbitrary dimension. More precisely, given a morphism of pairs
\[
\sigma : (X,E) \to (Y,F),
\]
where $X$ and $Y$ are smooth varieties over an algebraically closed field of characteristic zero, and $E$ and $F$ are simple normal crossings divisors satisfying $\sigma^{-1}(F)_{\mathrm{red}} \subset E$. We write the dimension of $X$ as $n$ and that of $Y$ as $N$. We associate the logarithmic Fitting ideals
\[
\cF_{n-k}(\sigma), \qquad k\in\{1,\ldots,\mbox{min}\{n,N\}\},
\]
obtained from the logarithmic Jacobian matrix of $\sigma$. The principalisation of these ideals is expected to govern the simplification of the pulled-back cotangent sheaf and, ultimately, the monomialisation of the morphism.

The main problem considered in this paper can now be formulated as follows.

\begin{problem}[Principalisation of the second log-Fitting ideal]\label{mainpro}
Let
\[
\sigma : (X,E) \to (Y,F)
\]
be a morphism of pairs. Find a family of local proper birational morphisms
\[
\tau : (X',E') \to (X,E),
\]
each given by a finite sequence of smooth local blow-ups compatible with the divisor $E$, such that the first two logarithmic Fitting ideals
\[
\mathcal F_{n-1}(\sigma \circ \tau)
\qquad \text{and} \qquad
\mathcal F_{n-2}(\sigma \circ \tau)
\]
become principal monomial ideals.
\end{problem}

Positive results are already known in low dimensions, for example when $\mbox{dim}(X) \leq 3$ and $\sigma$ is a resolution of singularities \cite{da2017resolution}, or when $\mbox{dim}(X)\leq 3$ and $\mbox{dim}(Y)\leq 2$ \cite{cutkosky2002monomialization}. Furthermore, for dominant morphisms from $3$-folds to surfaces, Cutkosky established the existence of monomialisation in \cite{cutkosky2002monomialization,Cutkosky2013,Cut24}, and  the relationship between strong preparation
and logarithmic Fitting ideals is clarified
in \cite[Proposition 6.6]{jiang2025newproofmonomialisation3folds};
see also \cite{cutkosky2002monomialization,cutkosky2004monomialization,bbmono}. These results strongly suggest that the second logarithmic Fitting ideal captures a decisive obstruction to monomiality.

The goal of the present paper is to provide a local
resolution theorem for the first two logarithmic
Fitting ideals in arbitrary dimension. Our approach
builds on ideas going back to Cutkosky's work on
the monomialisation of morphisms from 3-folds to
surfaces \cite{cutkosky2002monomialization}. It draws on the techniques
developed in \cite{da2018local} for the local
monomialisation of systems of first integrals
of Darboux type, as well as the framework of
logarithmic differential ideals and
Weierstrass--Tschirnhausen preparation developed
in \cite{bbmono}. Suitable local blow-ups, together with a preparation procedure
of Weierstrass--Tschirnhausen type, allow us to reduce the complexity
of the logarithmic Jacobian matrix inductively. The main theorem of
the present paper gives a positive local answer in arbitrary dimension;
see Theorem \ref{thm:Local}.

A central role in our construction is played by the notion of \emph{log-rank adapted morphism}. Roughly speaking, such morphisms have a sufficiently controlled logarithmic rank stratification to allow the construction of local invariants compatible with the logarithmic structure. This notion first appeared naturally in the author’s Ph.D. thesis \cite{jiang2025}, where it served as a framework for studying the resolution of cotangent sheaves via logarithmic Fitting ideals.

The main result of this article is the following local principalisation theorem.

\begin{theorem}[Main Theorem]\label{thm:main}
Let
\[
\sigma : (X,E) \to (Y,F)
\]
be a log-rank adapted morphism of pairs. Then there exists a finite family of commutative diagrams
\[
\begin{tikzcd}
(X_\lambda,E_\lambda) \arrow[r,"\tau_\lambda"] \arrow[dr,"\sigma_\lambda"'] & (X,E) \arrow[d,"\sigma"] \\
& (Y,F)
\end{tikzcd}
\]
such that:
\begin{enumerate}
    \item each $\tau_\lambda$ is a finite sequence of smooth local blow-ups compatible with $E$;
    \item the images of the morphisms $\tau_\lambda$ cover $X$;
    \item for every $\lambda$, the logarithmic Fitting ideals
    \[
    \mathcal F_{n-1}(\sigma_\lambda)
    \qquad \text{and} \qquad
    \mathcal F_{n-2}(\sigma_\lambda)
    \]
    are principal monomial ideals.
\end{enumerate}
\end{theorem}

The proof combines several ingredients. First, the log-rank adapted hypothesis allows us to develop a suitable local normal form for the morphism, analogous to the prepared forms appearing in the classical theory of resolution of singularities. More precisely, after admissible local coordinate changes, one obtains a Weierstrass--Tschirnhausen type presentation whose coefficients encode the failure of the logarithmic Fitting ideals to be principal monomial.

The strategy is then reminiscent of the inductive philosophy underlying Hironaka’s proof of resolution of singularities. One principalises coefficient ideals attached to the prepared form and studies how the associated invariants transform under local blow-ups. In this sense, the present work may be viewed as an incarnation of local resolution of singularities in the context of logarithmic cotangent sheaves.

Although the theorem proved here is local, it constitutes one of the fundamental steps towards a functorial global resolution algorithm for logarithmic Fitting ideals. Indeed, as already discussed in \cite{jiang2025}, the global problem is intimately related to the construction of canonical invariants stable under smooth morphisms. Such functoriality would in turn provide a systematic approach to monomialisation problems in higher dimensions.

The theorem proved here is local: it produces a covering by sequences
of admissible local blow-ups, rather than a single globally defined
sequence. Indeed, the preparation procedure uses local choices of
\'etale coordinates and hypersurfaces of maximal contact, and the
resulting centres are not shown to be compatible on overlaps. A
global, functorial principalisation would require an intrinsic way to
choose permissible centres and to control the corresponding invariants
under smooth base change. Weighted blow-ups, which have proved
effective in functorial embedded resolution \cite{ATW24}, may offer a
suitable framework for addressing these difficulties; see also
\cite{jiang2025}. Establishing such a global algorithm, however, is
beyond the scope of the present paper.

The paper is organised as follows. We recall some basic notions in Section~2. Section~3 introduces logarithmic Fitting ideals and log-rank. In Section~4, the class of log-rank adapted morphisms is defined. We then define the local invariant used throughout the paper. Finally, Section~5 presents the proof of the local principalisation theorem.

\section{Preliminaries}
\label{preli}
\subsection{Review of main definitions}

\begin{definition}[$\mathbb{K}$-variety]
\label{kv}
    Let $\mathbb{K}$ be an algebraically closed field of characteristic zero. A $\mathbb{K}$-variety is a separated, reduced scheme of finite type over $\mathbb{K}$, not necessarily irreducible.
\end{definition}


A $\mathbb{K}$-variety of dimension 1, 2 or 3 is called a curve, a surface or a threefold, respectively.

\medskip

Note that a subvariety of a $\mathbb{K}$-variety is still a $\mathbb{K}$-variety.

\begin{definition}
    A blow-up is the blow-up of a smooth subvariety of a $\mathbb{K}$-variety.
\end{definition}

\begin{remark}
    A blow-up is projective for $\mathbb{K}$-varieties. Hence, a $\mathbb{K}$-variety remains a $\mathbb{K}$-variety after blow-ups. Moreover, a projective morphism of Noetherian schemes is proper, hence of finite type.
\end{remark}

In Problem \ref{mainpro}, we seek to principalise log-Fitting ideals. We therefore need to
consider the notion of a locally principal ideal sheaf. It turns out that we only need to
require that the ideal sheaf is locally principal on all the closed points. Therefore, from now on, all points of a $\mathbb K$-variety are assumed to be closed unless stated otherwise.

\medskip

Let us also recall several notions from \cite{cutkosky2004monomialization,bbmono}.

\begin{definition}[Morphism of pairs]\label{def:MorphismPairs}
    Let $\sigma:X\rightarrow Y$ be a morphism of finite type from a smooth $\mathbb{K}$-variety $X$ to a smooth $\mathbb{K}$-variety $Y$ with reduced simple normal crossings (SNC) divisors $E$ on $X$ and $F$ on $Y$ such that $\sigma^{-1}(F)_{\text{red}}\subseteq E$. Then $\sigma$ is called a morphism of pairs, denoted by $\sigma:(X,E)\rightarrow (Y,F)$. We will write the dimension of $X$ as $n$ and the dimension of $Y$ as $N$.
\end{definition}

Let $X$ be a nonsingular $\mathbb{K}$-variety, $E$ be a simple normal crossings divisor on $X$. By \textbf{regular $E$-adapted coordinates} at a (closed) point $\pa\in X$, we mean a regular system of parameters
\[
\left(\bu,\bv\right) = (u_1,\ldots,u_s,v_1,\ldots,v_{n-s})
\]
in $\mathcal{O}_{X,\pa}$ such that $E$ is locally defined by $\{u_1\cdots u_s=0\}$ with $s\in\mathbb{N}^*$. In this case, we say that $\pa$ is an $s$-point. We will often use more general coordinates, namely \textbf{étale} or \textbf{formal}. If $h: (U,y)\rightarrow X$ is an étale neighbourhood of $\pa$, then $U$ is smooth
and $h^*(E)$ is a simple normal crossings divisor on $U$. If $(\bu,\bv)$ are regular $h^*(E)$-adapted coordinates in $\mathcal{O}_{U,y}$, then $(\bu,\bv)$ are called étale $E$-adapted coordinates at $\pa$. We also consider formal coordinates in the complete local ring
$\hat{\mathcal{O}}_{X,\pa}$.

Furthermore, an element
\[
f\in\mathbb{K}[u_1,\ldots,u_s,v_1,\ldots,v_{n-s}]
\]
is called \emph{monomial with respect to $E$} if $f\in \mathbb{K}[u_1,\ldots,u_s]$ and $f$ is a monomial. The same terminology applies when $f$ is an algebraic or formal power series.

On the other hand, we will denote $F$-adapted coordinates in the target $Y$ via the notation $\left(\bx,\by\right) =(x_1,\ldots,x_t,y_1,\ldots,y_{N-t})$, where $\{x_1 \cdots x_t=0\}$ is the local equation of the SNC divisor $F$ on $Y$.

An ideal sheaf $\mathcal I\subset\mathcal O_X$ is said to be
\emph{principal monomial with respect to $E$} if, at every point
$\pa\in X$, there exist $E$-adapted coordinates
$(u_1,\ldots,u_s,v_1,\ldots,v_{n-s})$ and a multi-index
$\bbe\in\mathbb N^s$ such that
\[
  \mathcal I\cdot\mathcal O_{X,a}=(u^{\bbe}),
  \qquad u^{\bbe}:=u_1^{\beta_1}\cdots u_s^{\beta_s}.
\]
In particular, the unit ideal is principal monomial.

Most of the time, we will omit the expression `with respect to $E$' when this is clear from the context.

We now introduce a special class of blow-ups that will be used later.

\begin{definition}[Combinatorial blow-up]\label{def:Combinatorial}
A combinatorial blow-up is a blow-up whose centre is an intersection subspace of the divisor $E$, which is everywhere locally given by $u_1=\cdots=u_k=0$ for $1\leq k\leq s$, with respect to $E$-adapted coordinates.
\end{definition}

We use the following notion of compatibility between a blow-up and a divisor.

Let \(E=\{u_{1}\cdots u_{s}=0\}\subset X\) be a simple normal crossings (SNC) divisor,
locally written in $E$-adapted coordinates \(\left(\bu,\bv\right)=(u_{1},\dots,u_{s},v_{1},\dots,v_{n-s})\).

\begin{definition}[Blowing-up compatible with a divisor]\label{Def:compatible}
A (local) blow-up
\[
\sigma:X' \longrightarrow X
\]
with centre \(C\subset X\) is \emph{compatible with the divisor \(E\subset X\)} if $\forall\pa\in E\cap C$, there are $E$-adapted coordinates at $\pa$ for which
\begin{enumerate}
\item each irreducible component of \(E\) is a coordinate hypersurface
\(\{u_i=0\}\);
\item the centre \(C\) is a coordinate subspace, i.e.\ it is given by
the simultaneous vanishing of some of the coordinates.
\end{enumerate}
Equivalently, the pair \((E,C)\) itself has simple normal crossings.
In this situation the reduced total transform
\[
E' \;=\;\sigma^{-1}(E)_{\text{red}}
\]
remains an SNC divisor on \(X'\).
\end{definition}


Compatibility ensures that every step of the resolution keeps the divisor in
normal-crossings form. Consequently,
\begin{itemize}
\item the new simple normal crossings divisor \(E'\) is still SNC after each blow-up;
\item centres and transforms can be described purely in terms of coordinate
hyperplanes, which is what ultimately allows one to ``monomialise’’ a
morphism in suitable coordinates.
\end{itemize}

\begin{example}
In \(\mathbb{R}^{3}\), let \(E=\{x\,y=0\}\), the union of the planes \(x=0\) and
\(y=0\).  By definition, blowing up the line \(C=\{x=y=0\}\) (the \(z\)-axis) is compatible.
Blowing up an arbitrary smooth curve that meets \(E\) non-transversely is
\emph{not} compatible, since no coordinate system makes that curve a coordinate
subspace at the crossing points.
\end{example}

\section{Log-Fitting ideals and log-rank}\label{ssec:LogFittingCharacterizationStronglyPrepared}
\subsection{Definition of log-Fitting ideals}
Let $\left(X,E\right)$ be a smooth $\mathbb{K}$-variety of dimension $n$ along with a reduced SNC divisor. We denote by $\Omega_X^1$ the $\mathcal{O}_X$-module of 1-forms on $X$ and by $\Omega_X^1\left(\text{log}E\right)$ the $\mathcal{O}_X$-module of logarithmic 1-forms on $X$. Let $\pa\in X$ and let $\left(\bu,\bv\right)$ be \'etale $E$-adapted coordinates at $\pa$. Thus, there exists an étale neighbourhood $h:(U,y)\rightarrow X$ of $\pa$ such that $(\bu,\bv)$ are regular $h^*(E)$-adapted coordinates in $\mathcal{O}_{U,y}$. Then $h^*\Omega^1_X(logE)=\Omega^1_U(logh^*(E))$ is a locally
free $\mathcal{O}_{U,y}$-module, and after replacing $U$ with a Zariski open neighbourhood of $y$ in
$U$, $\Omega^1_U(logh^*(E))$ has the $\mathcal{O}_U$-basis
\[
\frac{du_i}{u_i},i=1,\ldots,s,\text{ and }dv_j,j=1,\ldots,n-s.
\]
Let $\sigma:(X,E)\rightarrow(Y,\emptyset)$ be a morphism of pairs with dim$Y=N$, and let $\sigma$ be expressed as $(\sigma_1,\ldots,\sigma_N)$. The logarithmic Jacobian matrix of $\sigma$, analogous to the Jacobian matrix of $\sigma$, is defined using the derivations dual to the basis consisting of $\frac{du_i}{u_i}$ and $dv_j$. More precisely,
\[
\mbox{logJac}\sigma=
\begin{bmatrix}
u_1\partial_{u_1}\sigma_1 & \ldots & u_s\partial_{u_s}\sigma_1 & \partial_{v_1}\sigma_1 & \ldots & \partial_{v_{n-s}}\sigma_1\\
u_1\partial_{u_1}\sigma_2 & \ldots & u_s\partial_{u_s}\sigma_2 & \partial_{v_1}\sigma_2 & \ldots & \partial_{v_{n-s}}\sigma_2\\
\vdots & & \vdots & \vdots & & \vdots\\
u_1\partial_{u_1}\sigma_N & \ldots & u_s\partial_{u_s}\sigma_N & \partial_{v_1}\sigma_N & \ldots & \partial_{v_{n-s}}\sigma_N\\
\end{bmatrix}
\]

The sheaf $\Omega^k_X\left(\text{log}E\right)$ of logarithmic $k$-forms over $X$ is defined in terms of $\Omega^1_X\left(\text{log}E\right)$ in the standard way.

Let $\sigma:(X,E)\rightarrow (Y,F)$ be a morphism of pairs. For all $k=1, \ldots, N$, $\sigma^*\left(\Omega_Y^k(\log F)\right)$ is a subsheaf of $\Omega_X^k(\log E)$. For the proof, see {\cite[Lemma 4.1]{bbmono}}.

We now define the logarithmic Fitting ideals, which play a central role in this paper.

\begin{definition}[Logarithmic Fitting ideals of order $k$]
For any $k\in\{1,\ldots,\mbox{min}\{n,N\}\}$, the logarithmic Fitting ideal of order $k$ is the subsheaf of $\mathcal{O}_X$ locally generated by the minors of order $k$ of logJac$\sigma$ with respect to regular $E$-adapted coordinates at $\pa$, and we denote it by $\mathcal{F}_{n-k}(\sigma)$. We shall abbreviate logarithmic Fitting ideals as log-Fitting ideals.
\end{definition}

If $h:(U,y)\rightarrow X$ is an étale neighbourhood of $\pa$, then
\[
\mathcal{F}_{n-k}(\sigma\circ h)=\mathcal{F}_{n-k}(\sigma)\cdot\mathcal{O}_U
\]
and logrk$_y(\sigma\circ h)=\mbox{logrk}_{\pa}(\sigma)$.

\subsection{Definition of log-rank}
\begin{definition}\label{def:LogRank}
Let $\sigma:(X,E)\rightarrow (Y,F)$ be a morphism of pairs. Given $\pa\in X$, the logarithmic rank of $\sigma$ at $\pa$ is defined as 
\[
\mbox{logrk}_{\pa}(\sigma) := \sup \{k ;\, \cF_{n-k}(\sigma) \cdot \mathcal{O}_{X,\pa} = \mathcal{O}_{X,\pa}\},
\]
and is defined to be zero if the above set is empty.
\end{definition}

Equivalently, the log-rank can be characterised as follows.

\begin{definition}
    Let $\sigma:(X,E)\rightarrow (Y,F)$ be a morphism of pairs. $\forall\pa\in X$, the logarithmic rank of $\sigma$ at $\pa$ is defined as the rank of the logarithmic Jacobian matrix evaluated at $\pa$, and we shall denote it by \mbox{logrk}$_{\pa}\left(\sigma\right)$. If $h:(U,y)\rightarrow X$ is an étale neighbourhood of $\pa$, then logrk$_y(\sigma\circ h)=\mbox{logrk}_{\pa}(\sigma)$.
\end{definition}

\section{Log-rank adapted morphism and local invariant}
\label{lam}
We introduce the notion of a log-rank adapted morphism, following discussions with A. Belotto. This notion is essential to our proof and is also useful in other contexts.

Recall the definition of quasi-prepared morphism.

\begin{definition}[Quasi-prepared morphism, see {\cite[Definition 3.1]{cutkosky2004monomialization}}]\label{def:QuasiPrepared}
    Let $\sigma:(X,E)\rightarrow (Y,F)$ be a morphism of pairs (not necessarily dominant), and let
    \[
 \operatorname{Sing}(\sigma)
 :=\{\pa\in X: \operatorname{rank}(d\sigma_{\pa})<\min\{n,N\}\},
\]
be the singular locus of the morphism. Here $d\sigma_{\pa}:T_{\pa}X\to T_{\sigma(\pa)}Y$ is the ordinary
differential of $\sigma$, and the displayed set is the singular locus of $\sigma$. We say that $\sigma$ is quasi-prepared if Sing$\left(\sigma\right)\subseteq E$ and $\sigma^{-1}(F)_{\text{red}}=E$. 
\end{definition}

Let $\sigma:(X,E)\rightarrow (Y,F)$ be a morphism of pairs. For any $i\in 1,\ldots,\mbox{min}\{n,N\}$, note $\Sigma_i=\{\pa\in E:\mbox{logrk}_{\pa}\left(\sigma\right)\leq \mbox{min}\{n,N\}-i\}$ and $Y_i=\overline{\sigma\left(\Sigma_i\right)}$. Observe that $\Sigma_i$ is a $\mathbb{K}$-variety by definition. The ideal sheaf of $Y_i$ is denoted by $\mathcal{I}_{Y_i}$.

\begin{definition}[Log-rank adapted morphism]\label{wm}
Let $\sigma:(X,E)\rightarrow (Y,F)$ be a quasi-prepared morphism (not necessarily dominant), where $\mbox{dim}X=n$ and $\mbox{dim}Y=N$. Denote by $\sigma_0:(X,E) \rightarrow (Y,\emptyset)$ the morphism where the divisor $F$ is declared to be empty, and we shall use $\sigma_0$ to define $\Sigma_i$ as well as $Y_i$.

We say that $\sigma$ is log-rank adapted at an $s$-point $\pa\in E$ if we can choose regular $E$-adapted coordinates $\left(\bu,\bv\right)$ at $\pa$ and regular $F$-adapted coordinates $(\bw,\bx)$ at $\sigma\left(\pa\right)$ such that
\begin{enumerate}
\item $\sigma=\left(\sigma_1,\ldots,\sigma_N\right)$ satisfies
$$
\sigma_1=v_1,\ldots,\sigma_r=v_r,\sigma_{r+1}=\varepsilon\bu^{\bal},
$$
where $r\in\mathbb{N}$ is the log-rank of $\sigma_0$ at $\pa$, $\varepsilon\in\mathcal O_{X,\pa}^\times$, $\bal\in\mathbb{N}^s$, $\mathcal{I}_{Y_{N-r}}$ is generated by $z_{r+1},\ldots,z_N$ at $\sigma\left(\pa\right)$ and $\sigma^*\mathcal{I}_{Y_{N-r}}$ is generated by $\sigma_{r+1}=\sigma^*\left(z_{r+1}\right)$ at $\pa$.
\item There is a finite sequence of SNC divisors $E^{(k)}$ such that $E^{(1)}=E$, $E^{(k+1)}\subseteq E^{(k)}$, $\forall k=1,\ldots,\mbox{min}\{n,N\}$ and
\begin{enumerate}
\item logrk$_{\pa}\left(\sigma\right)=\mbox{min}\{n,N\}-k$, $\forall \pa\in E^{(k)}\setminus E^{(k+1)}$;
\item $\overline{\sigma\left(E^{(k)}\right)}\setminus\overline{\sigma\left(E^{(k+1)}\right)}$ is smooth;
\item
$Y_i$ has simple normal crossings with respect to $F$, $\forall k=1,\ldots,\mbox{min}\{n,N\}$.
\end{enumerate}
\end{enumerate}
\end{definition}

For a log-rank adapted morphism $\sigma:(X,E)\rightarrow (Y,F)$, we may further choose étale $E$-adapted coordinates $(\bu,\bv)$ at a point $\pa\in X$ such that $\sigma_{r+1}=\bu^{\bal}$.

The relevance of this class is demonstrated by the reduction theorem
of \cite[Theorem~0.3.1]{jiang2025}:
\begin{theorem}[Reduction to log-rank adapted morphism]\label{thm:LrAdapted}
Let $\sigma: (X,E) \to (Y,F)$ be a morphism of pairs. Suppose that one of the following conditions holds:
\begin{itemize}
\item $\sigma$ is a dominant morphism such that \(
\overline{\sigma\left(E\cup\text{Sing}(X)\right)}
\) has codimension $\geq 1$.

\item $\sigma$ is a quasi-prepared morphism, see Definition \ref{def:QuasiPrepared};
\end{itemize}
Then there exist finite sequences of blow-ups $\tau^{\prime}: (X',E')\to (X,E), \tau: (Y',F') \to (Y,F)$, and a log-rank adapted morphism $\sigma':(X',E') \to (Y',F')$, see Definition \ref{wm}, such that the following diagram commutes:
\begin{center}
\begin{tikzcd}
X^{\prime} \arrow[r, "\sigma^{\prime}"] \arrow[d, "\tau^{\prime}"]
& Y^{\prime} \arrow[d, "\tau"] \\
X \arrow[r, "\sigma"]
& Y.
\end{tikzcd}
\end{center}
This reduction is functorial with respect to smooth morphisms to the target.
\end{theorem}

In other words, under the hypotheses stated above,
a morphism of pairs can be transformed into a log-rank adapted
morphism by finite sequences of blow-ups of the source
and the target.

\begin{definition}[Log-derivation with respect to an exceptional divisor]
Let $\mathcal D_E$ be the sheaf of logarithmic $\mathbb K$-derivations
along $E$ on a $\mathbb K$-variety $X$. In $E$-adapted coordinates, it is locally generated
by
\[
 u_1\partial_{u_1},\ldots,u_s\partial_{u_s},
 \partial_{v_1},\ldots,\partial_{v_{n-s}}.
\]
For a coherent ideal sheaf $\mathcal I\subset\mathcal O_X$
and $k\ge 0$, let $\mathcal D_E^{\le k}(\mathcal I)$ be the
ideal locally generated by all expressions
\[
 \delta_1\cdots\delta_j(f),\qquad
 0\le j\le k,\quad f\in\mathcal I,\quad
 \delta_i\in\mathcal D_E.
\]
We set
\[
 \mathcal D_E^{\le\infty}(\mathcal I)
 :=\sum_{k\ge0}\mathcal D_E^{\le k}(\mathcal I).
\]
This is the smallest $\mathcal D_E$-invariant ideal containing
$\mathcal I$; it is locally finitely generated because $X$ is
Noetherian.
\end{definition}

\begin{remark}
By the definition of a log-rank adapted morphism, at any point $\pa$ where $\operatorname{logrk}_{\pa}(\sigma_0)\ge1$,
the second log-Fitting ideal is already principal monomial.
Indeed, if $\operatorname{logrk}_{\pa}(\sigma_0)\ge2$, then
$\mathcal F_{n-2}(\sigma)\cdot\mathcal O_{X,\pa}=\mathcal O_{X,\pa}$; if $\operatorname{logrk}_{\pa}(\sigma_0)=1$, one may choose étale $E$-adapted coordinates such that
$\sigma_1=v_1$, $\sigma_2=\bu^{\bal}$, and
$\sigma_i\in(\bu^{\bal})$ for $i\ge2$.
Since the logarithmic derivatives of the latter functions are
divisible by $\bu^{\bal}$, all logarithmic $2\times2$ minors
are divisible by $\bu^{\bal}$. Conversely, for any $j$
with $\alpha_j>0$, the minor involving the rows
$\sigma_1,\sigma_2$ and the columns
$\partial_{v_1},u_j\partial_{u_j}$ equals
$\pm\alpha_j \bu^{\bal}$. Consequently,
\[
 \mathcal F_{n-2}(\sigma)\cdot\mathcal O_{X,\pa}=(\bu^{\bal}).
\]
Thus, the nontrivial part of the principalisation problem
is concentrated on the log-rank zero locus.
\end{remark}

Now, given a log-rank adapted morphism $\sigma_0: (X_0,E_0) \to (Y,F)$, consider the following ideal sheaf:
\[
\mathcal{J}_0:=D_{E_0}^{\leq +\infty}\left(\mathcal{F}_{n-2}\left(\sigma_0\right)\right),
\]
and it is straightforward to check that $\mathcal{J}_0$ is generated by monomials with support in $E_0$.

\begin{lemma}\label{Lem:PJ}
There exists a sequence of combinatorial blow-ups $\tau_0: (X,E) \to (X_0,E_0)$ such that
\[
\mathcal{J} =D_E^{\leq +\infty}\left(\mathcal{F}_{n-2}\left(\sigma_0 \circ \tau_0\right)\right)
\]
is principal monomial.
\end{lemma}
\begin{proof}
We seek to resolve $\mathcal{J}_0$. Since $\mathcal{J}_0$ is closed under logarithmic derivations, $D_{E_0}(\mathcal{J}_0)\subset\mathcal{J}_0$. Therefore, all the centres to resolve $\mathcal{J}_0$ are $D_{E_0}$-invariant by a result of A. Belotto, see \cite[Theorem 1.4.2]{da2013resolution}. In other words, if $\mathcal{I}$ is the ideal of a centre used to resolve $\mathcal{J}_0$, it satisfies $D_{E_0}(\mathcal{I})\subset\mathcal{I}$. In particular, this implies that $v_j\notin\mathcal{I}$ for $j=1,\ldots,n-s$, since otherwise $D_{E_0}(\mathcal{I})=\mathcal{O}_X$. Then all the blow-ups in the resolution of $\mathcal{J}_0$ are combinatorial.
\end{proof}

\begin{definition}[Order with respect to an exceptional divisor] Let $\mathcal{I}$ be an ideal sheaf. We define:
\[
\mbox{ord}_{E,\pa}(\mathcal{I}) = \mbox{inf} \, \{k \in \mathbb{N}\mbox{ s.t. }\mathcal{O}_{X,\pa}\cdot D^{\leq k}_{E,\pa}\left(\mathcal{I}\right)=\mathcal{O}_{X,\pa}\}.
\]
For a local section $f\in\mathcal{O}_{X,\pa}$, the order of $f$ with respect to the exceptional divisor $E$ is defined as $\mbox{ord}_{E,\pa}(\mathcal{I})$ for $\mathcal{I}=(f)$, and is denoted by $\mbox{ord}_{E,\pa}(f)$.
\end{definition}

Lemma \ref{Lem:PJ} allows us to obtain a more convenient local presentation of $\sigma$:


\begin{lemma}[Weierstrass-Tschirnhausen normal form, see {\cite[Definition~4.3.1, Theorem~4.5.1]{jiang2025}}; see also {\cite[Lemmas 5.23 and 5.31]{bbmono}} for related Weierstrass–Tschirnhausen normal forms arising from logarithmic differential closures.]\label{Lemma:sigmalogWT}
Let $\sigma: (X,E) \to (Y,F)$ be a log-rank adapted morphism, and let $\pa \in X$ be a point at which the log-rank is $0$. If $\mathcal{J}$ is principal monomial at $\pa$ and $\mathcal{F}_{n-2} \cdot \mathcal{O}_{X,\pa}\neq \mathcal{J}\cdot \mathcal{O}_{X,\pa}$, then there exist étale $E$-adapted coordinates $\left(\bu,\bv\right)$ in which $\sigma$ has \textbf{Weierstrass-Tschirnhausen normal form}, i.e.
$$
\left\{
\begin{aligned}
\sigma_1&=\bu^{\bal}, \bal\in \mathbb{N}^r \\
\sigma_i&=g_i\left(\bu\right)+\bu^{\bde}T_i, i=2,\ldots,N \\
\end{aligned}
\right.
$$
where
\begin{enumerate}
\item $\forall i$, $\bu^{\bal}$ divides $\sigma_i$, and the ideal $\mathcal T_{\left(\bu,\bv\right)}:=(T_2,\ldots,T_N)$ (which depends on the coordinates $\left(\bu,\bv\right)$) cannot be factored by $\bu$;
\item $g_i\left(\bu\right)\in \mathbb{K}[|\bu|]$ such that $d\bu^{\alpha}\wedge dg_i\left(\bu\right)=0$, for $i=2,\ldots,N$;
\item All exponents of the monomials in the formal expansion of $\bu^{\bde}T_i\in \mathcal{O}^{h}_{X,\pa}[|\bu|]$ are linearly independent with $\bu^{\bal}$, for $i=2,\ldots,N$.
\end{enumerate}

Moreover,
$$ \left\{
\begin{aligned}
T_2 & =\tilde{T}_2\left(\bu,\bv\right)v_1^{d_1}+\sum_{j=0}^{d_1-2}a_{2j}\left(\bu,\hat{\bv}\right)v_1^{j}, \\
T_i & = \tilde{T}_i\left(\bu,\bv\right)v_1^{d_1}+\sum_{j=0}^{d_1-1}a_{ij}\left(\bu,\hat{\bv}\right)v_1^j,i\in \{3,\ldots,N\}, \\
\end{aligned}
\right.
$$
where $\tilde{T}_2$ is a unit, $a_{ij}\left(\bze\right)=0$, $\forall i\in\{2,\ldots,N\},j\in\{0,\ldots,d_1-1\}$, $\hat{\bv}$ denotes the tuple obtained from $\bv$ by omitting $v_1$, and:
\[ d_1(\bu,\bv):=\mbox{min}\{\mbox{ord}_{E,\pa}\left(T_i(\bu,\bv)\right)\}_{i=2,\ldots,N} >1.
\]
\end{lemma}

At a log-rank zero point $\pa$, if
$\mathcal F_{n-2}(\sigma)\cdot\mathcal O_{X,\pa}
=\mathcal J\cdot\mathcal O_{X,\pa}\ne\mathcal O_{X,\pa}$,
we set $d_1=1$ at $\pa$.

If
$\mathcal F_{n-2}(\sigma)\cdot\mathcal O_{X,\pa}
=\mathcal O_{X,\pa}$,
we set $d_1=0$ at $\pa$.

Otherwise, we define
\[
d_1=\min_{2\le i\le N}
\operatorname{ord}_{E,a}(T_i)
\]
using the Weierstrass--Tschirnhausen
presentation of Lemma \ref{Lemma:sigmalogWT}.

\section{Local resolution}
We provide a solution to Problem \ref{mainpro} via \emph{local blow-ups}, following the methods of \cite{bbmono,da2018local,cutkosky2002monomialization}.

Recall that a smooth blow-up is a blow-up with smooth centre.

\begin{definition}[Local blow-up, see {\cite[Definitions and Remarks 1.9]{bbmono}}]\label{Def:LocBlow}
In the complex- or real-analytic setting, a (smooth) local blow-up $\sigma: X^{\prime} \rightarrow X$ is the composition $\sigma_1 \circ \sigma_2$ of an open immersion $\sigma_1: U \hookrightarrow X$ and a (smooth) blow-up $\sigma_2: X^{\prime} \rightarrow U$.

In the algebraic case, an open subset is understood as an étale open, so that a local blow-up $\sigma$ means the composite $\sigma_1 \circ \sigma_2$ of an étale morphism $\sigma_1: U \rightarrow X$ and a smooth blow-up $\sigma_2: X^{\prime} \rightarrow U$.
\end{definition}

In this section, we write $(\bu,\bv)=(u_1,\ldots,u_r,v_1,\ldots,v_s)$, $r+s=n$ for \'{e}tale $E$-adapted coordinates; the indices have been relabelled for convenience.

Our main theorem is:
\begin{theorem}\label{thm:Local}
Let $\sigma:(X,E)\rightarrow (Y,F)$ be an algebraic log-rank adapted morphism. There exists a finite family of diagrams
\[
\begin{tikzcd}[column sep=large, row sep=huge]
\left(X_{\lambda},E_{\lambda}\right) \arrow[d, "\tau_{\lambda}"'] \arrow[rd, "\sigma_{\lambda}"] & \\ (X,E) \arrow[r, "\sigma"'] & (Y,F)
\end{tikzcd}
\]
such that
\begin{enumerate}
\item each $\tau_{\lambda}$ is a finite sequence of smooth local blow-ups, compatible with $E$;
\item for each $\lambda$, denote the image of $\tau_{\lambda}$ by Im$\left(\tau_{\lambda}\right)$; then $\cup_{\lambda}\mbox{Im}\left(\tau_{\lambda}\right)=X$;
\item the log-Fitting ideals $\mathcal{F}_{n-1}\left(\sigma_{\lambda}\right),\mathcal{F}_{n-2}\left(\sigma_{\lambda}\right)$ are principal monomial.
\end{enumerate}
\end{theorem}

\begin{remark}
If $\sigma:(X,E)\rightarrow (Y,F)$ is an analytic log-rank adapted morphism, then the family $\tau_{\lambda}$ is not finite, but it is \emph{semi-proper}, that is, for every relatively compact set $K \subset X$, there exist finitely many $\tau_{\lambda}$ and corresponding compact sets $K_{\lambda}\subset X_{\lambda}$ such that the union of $\tau_{\lambda}(K_{\lambda})$ covers $K$. We give the proof in the algebraic setting; it can be adapted directly to the analytic setting. 
\end{remark}


By Lemma \ref{Lem:PJ}, without loss of generality, one can assume that $\mathcal{J}_0$ is principal monomial, so $d_1$ is well-defined.

\begin{theorem}\label{Theo:main}
Let $(X,E)$ be a pair where $X$ is a $\mathbb{K}$-variety, a real-analytic or a complex-analytic space and $E$ a SNC divisor on $X$. For any point $\pa\in X$, there exists a family of $\tau_{\lambda}:\left(X_{\lambda},E_{\lambda}\right)\rightarrow (X,E)$ satisfying:
\begin{enumerate}
\item each $\tau_{\lambda}$ is a finite sequence of smooth local blow-ups, compatible with $E$;
\item the union $\cup_{\lambda}\mbox{Im}\left(\tau_{\lambda}\right)$ covers a neighbourhood of $\pa$;
\item at any point $\pb\in\tau_{\lambda}^{-1}\left(\pa\right)$, the invariant $d_1$ is 0 or 1.
\end{enumerate}
\end{theorem}
\begin{proof}
We argue by induction on the dimension of $X$.

When dim$X$=0, $X$ is a finite union of isolated points, so the assertion is immediate.

Assume that the assertion is true for any dimension strictly less than $\mathbb{D}$ with $\mathbb{D}\in\mathbb{N}^*$. It remains to prove Propositions \ref{Prop:pre} and \ref{Prop:de}.
\end{proof}

\begin{proposition}[Preparation]\label{Prop:pre}
Let $(X,E)$ be a pair where $X$ is a $\mathbb{K}$-variety, a real-analytic or a complex-analytic space of dimension $\mathbb{D}\in\mathbb{N}^*$ and $E$ a SNC divisor on $X$. For any point $\pa\in X$, suppose that the invariant $d_1$ at $\pa$ is finite, $d_1>1$ and Theorem \ref{Theo:main} is valid for any dimension strictly less than $\mathbb{D}$. Then there exists a family of $\tau_{\lambda}:\left(X_{\lambda},E_{\lambda}\right)\rightarrow (X,E)$ satisfying:
\begin{enumerate}
\item each $\tau_{\lambda}$ is a finite sequence of smooth local blow-ups, compatible with $E$;
\item the union $\cup_{\lambda}\mbox{Im}\left(\tau_{\lambda}\right)$ covers a neighbourhood of $\pa$;
\end{enumerate}
such that, at any point $\pb\in\tau_{\lambda}^{-1}\left(\pa\right)$, there are étale $E$-adapted coordinates $\left(\bu,v_1,\hat{\bv}\right)$ and
$$
\left\{
\begin{aligned}
T_2&=\tilde{T}_2\left(\bu,\bv\right)v_1^{d_1}+\sum_{j=0}^{d_1-2}a_{2j}\left(\bu,\hat{\bv}\right)v_1^{j}=v_1^{d_1}U\left(\bu,v_1,\hat{\bv}\right)+\sum_{j=1}^{d_1-2}\bu^{\bga_{2j}}b_{2j}\left(\bu,\hat{\bv}\right)v_1^j+a_{20}\left(\bu,\hat{\bv}\right), \\
T_i & = \tilde{T}_i\left(\bu,\bv\right)v_1^{d_1}+\sum_{j=0}^{d_1-1}a_{ij}\left(\bu,\hat{\bv}\right)v_1^{j}=\tilde{T}_i\left(\bu,\bv\right)v_1^{d_1}+\sum_{j=0}^{d_1-1}\bu^{\bga_{2j}}b_{ij}\left(\bu,\hat{\bv}\right)v_1^j,i\in \{3,\ldots,N\},
\end{aligned}
\right.
$$
where $\tilde{T}_2\left(\bu,\bv\right)$ is a unit and either $b_{ij}\left(\bze\right)\neq 0$ or $b_{ij}\equiv 0$ for $i=2,j=1,\ldots,d_1-2$ and for $i=3,\ldots,N,j=1,\ldots,d_1-1$. Furthermore, there exists $i_0\in\{2,\ldots,N\}$ such that either $a_{i_00}\equiv 0$ or $a_{i_00}=\bu^{\bbe}v_s^{\epsilon}$ with $\epsilon\in\{0,1\}$ and $\bu^{\bbe}$ divides $a_{i0}$, $\forall i\in\{2,\ldots,N\}$.
\end{proposition}
\begin{proof}
Since the morphism $\sigma$ in Theorem \ref{thm:Local} is log-rank adapted, $\mathcal{F}_{n-1}(\sigma)$ is principal monomial. It therefore suffices to preserve these two properties under the following blow-ups. We proceed by induction on $d_1$ to make $\mathcal{F}_{n-2}(\sigma)$ principal monomial while preserving the two properties of $\mathcal{F}_{n-1}(\sigma)$.

Define $\Lambda:=\{(i,j);a_{ij}\not\equiv 0\}$ and $\mathcal{I}:=\prod_{(i,j)\in\Lambda}a_{ij}\left(\bu,\hat{\bv}\right)$. After resolving this ideal over $U$, we obtain the following morphism:
\[
\pi:\left(\tilde{U},\widetilde{E\cap U}\right)\rightarrow (U,E\cap U)
\]
such that $\pi^*\left(\mathcal{I}\right)$ is principal monomial. Consequently, $\forall\tilde{\pa}\in\pi^{-1}\left(\pa\right)$, we have
$$ \left\{
\begin{aligned}
\pi^*\sigma_1&=\tilde{\bu}^{\tilde{\bal}}, \tilde{\bal}\in \mathbb{N}^r \\
\pi^*\sigma_i&=\tilde{g}_i\left(\tilde{\bu}\right)+\tilde{\bu}^{\tilde{\bde}}T_i, i=2,\ldots,N \\
\end{aligned}
\right.
$$
and $\pi^*a_{ij}=\tilde{\bu}^{\tilde{\bga}_{ij}}b_{ij}$ with either $b_{ij}\left(\bze\right)\neq 0$ or $b_{ij}\equiv 0$, for $i\in\{2,\ldots,N\},j\in\{0,\ldots,d_1-1\}$.

Moreover, $\mathcal{F}_{n-1}(\pi\circ\sigma)$ is principal monomial because the blow-ups of $\pi$ are compatible with $E$ (Definition \ref{Def:compatible}).

After applying the resolution of singularities to $\mathcal{I}:=\prod_{(i,j)\in\Lambda}a_{ij}\left(\bu,\hat{\bv}\right)$ as described above, for any point $\pa\in X$, we use Weierstrass-Tschirnhausen normal form (defined in Lemma \ref{Lemma:sigmalogWT}) in an étale neighbourhood $U$ of $\pa$ in such a way that the étale coordinates of the Weierstrass-Tschirnhausen normal form are induced from a regular system of parameters at $\pa$. We have
$$
\left\{
\begin{aligned}
T_2&=\tilde{T}_2\left(\bu,\bv\right)v_1^{d_1}+\sum_{j=0}^{d_1-2}a_{2j}\left(\bu,\hat{\bv}\right)v_1^{j}=\tilde{T}_2\left(\bu,\bv\right)v_1^{d_1}+\sum_{j=1}^{d_1-2}\bu^{\bga_{2j}}b_{2j}\left(\bu,\hat{\bv}\right)v_1^j+a_{20}\left(\bu,\hat{\bv}\right), \\
T_i & = \tilde{T}_i\left(\bu,\bv\right)v_1^{d_1}+\sum_{j=0}^{d_1-1}a_{ij}\left(\bu,\hat{\bv}\right)v_1^{j}=\tilde{T}_i\left(\bu,\bv\right)v_1^{d_1}+\sum_{j=0}^{d_1-1}\bu^{\bga_{2j}}b_{ij}\left(\bu,\hat{\bv}\right)v_1^j,i\in \{3,\ldots,N\},
\end{aligned}
\right.
$$
where $\tilde{T}_2\left(\bu,\bv\right)$ is a unit because it is the pullback of a unit, and $a_{ij}\left(\bze\right)=0$, $\forall i\in\{2,\ldots,N\},j\in\{0,\ldots,d_1-1\}$ because they are pullbacks of non-units.

Consider $U|_{v_1=0}$, the restriction of $U$ to the hypersurface of maximal contact $v_1=0$. Applying the inductive hypothesis that Theorem \ref{Theo:main} holds in every dimension strictly less than $\mathbb{D}$, we obtain $\tau_{\lambda}:\left(M_{\lambda},F_{\lambda}\right)\rightarrow\left(U|_{v_1=0},U\cap E|_{v_1=0}\right)$. Then the invariant at every point of $\tau_{\lambda}^{-1}\left(\pa\right)$ is 0 or 1. Consequently, for each $\lambda$, there exists $i_0\in\{2,\ldots,N\}$ such that $\tau_{\lambda}^*a_{i_00}=\bu^{\bbe}v_s^{\epsilon}$ with $\epsilon\in\{0,1\}$ if $\tau_{\lambda}^*a_{i_00}$ is not zero. Meanwhile, $\bu^{\bga_{\bga_{ij}}}$ is still principal monomial because the blow-ups of $\tau_{\lambda}$ are compatible with $E$ (Definition \ref{Def:compatible}), and we still use $\bu^{\bga_{\bga_ij}}$ to denote them with a slight abuse of notation.

We can now extend $\tau_{\lambda}$ to a blow-up of $X$ by taking the
product of the centres of $\tau_{\lambda}$ with the $v_1$-direction, and the result follows. Note that $\bbe$ is intrinsic by the definition of Weierstrass-Tschirnhausen normal form throughout the resolution process, so $\bu^{\bbe}$ divides $a_{i0}$, $\forall i\in\{2,\ldots,N\}$.
\end{proof}


\begin{proposition}[The invariant $d_1$ decreases strictly]\label{Prop:de}
Use the setting and notation of Proposition \ref{Prop:pre}. Let $\pa\in X$ be a point satisfying the conclusion of Proposition \ref{Prop:pre}. Then there exists a finite sequence $\sigma_{\pa}$ of smooth local blow-ups (i.e. defined in an étale neighbourhood $U$ of $\pa$) compatible with $E$ such that, at any point $\pb\in\sigma_{\pa}^{-1}\left(\pa\right)$, the invariant $d_1$ decreases strictly.
\end{proposition}
\begin{proof}
Recall our notation of étale $E$-adapted coordinates at $\pa$:

\(\left(\bu,\bv\right)=(u_{1},\dots,u_{r},v_{1},\dots,v_{s})\), where $r+s=n$. We now adopt the notation of {\cite[Section 10.1]{da2018local}} by setting
$$
\begin{aligned}
u_i & =u_i,i=1,\ldots,r \\
v & =v_1 \\
w_{r+j} & =v_j,j=2,\ldots,s
\end{aligned}
$$

We thus use $(\mathbf{u}, v, \mathbf{w})$ as our new notation.

We identify the dimension of $X$, denoted by $n$, with $m$ in the notation of {\cite[Section 10.1]{da2018local}}.

We use this notation for the remainder of the proof. After a sequence of local combinatorial blow-ups at $\pa$ with respect to $\left(u_1\ldots u_rv_1=0\right)$, we obtain new local coordinates at every point $\pb$ in the inverse image of $\pa$, which we denote by $\left(\bx,\bw\right)$. By definition, we have:
$$
(\mathbf{u}, v, \mathbf{w})=\left(\mathbf{x}^{\mathcal{A}}, \mathbf{w}\right)=\left(\mathbf{x}^{\mathbf{A}}, \mathbf{x}^\alpha, \mathbf{w}\right),
$$

Explicitly,
$$
\begin{aligned}
u_i & =x_1^{a_{j, 1}} \cdots x_{r+1}^{a_{j, r+1}},i=1,\ldots,r \\
v & =x_1^{\alpha_1} \cdots x_{r+1}^{\alpha_{+1}} \\
w_j & =w_j,j=r+2,\ldots,m
\end{aligned}
$$
with $\alpha_{i, j} \in \mathbb{N}$. Hence, the matrix $\mathcal{A}:=\left[\begin{array}{l}\mathbf{A} \\ \boldsymbol{\alpha}\end{array}\right]$ is a square matrix of size $(r+1)$ given by
$$
\mathbf{A}=\left[\begin{array}{ccc}
a_{1,1} & \ldots & a_{1, r+1} \\
\vdots & \ddots & \vdots \\
a_{r, 1} & \ldots & a_{r, r+1}
\end{array}\right] \quad \text { and } \quad \boldsymbol{\alpha}=\left[a_{r+1,1}, \ldots, a_{r+1, r+1}\right].
$$

Without loss of generality, we may translate the coordinates by setting $y_i=x_i-\xi_i$, $\xi_i\neq0$ for $i=t,\ldots,r+1$. Note that $t \neq 0$ as $\pb$ is in the pre-image of $\pa$. Hence, we obtain the coordinate system $(\mathbf{x}, \mathbf{y}, \mathbf{w})=\left(x_1, \ldots, x_t, y_{y+1}, \ldots, y_{r+1}, w_{r+2}, \ldots\right.$, $\left.w_m\right)$ centred at $\pb$.

We have a decomposition of $\mathcal{A}$:
$$
\mathcal{A}=\left[\begin{array}{ll}
\mathbf{A}_1 & \mathbf{A}_2 \\
\boldsymbol{\alpha}_1 & \boldsymbol{\alpha}_2
\end{array}\right]
$$
where $\mathbf{A}_1$ is an $r \times t$ matrix, $\mathbf{A}_2$ is an $r \times(r+1-t)$ matrix, $\boldsymbol{\alpha}_1$ is a $1 \times t$ matrix and $\boldsymbol{\alpha}_2$ is a $1 \times(r-t+1)$ matrix. We remark that, since $\pb$ is a point on the exceptional divisor, there exists at least one $u_i$ such that the pullback of $u_i(q)$ vanishes, which implies that $\mathbf{A}_1$ is a non-zero matrix. We distinguish two cases according to the rank of $\mathbf{A}_1$.

\begin{lemma}[Case 1, see {\cite[Lemma 10.2]{da2018local}}]\label{Lem:c1}
Assume $\mathbf{A}_1$ has maximal rank. Then, there exists a monomial coordinate system $(\mathbf{x}, \mathbf{y}, z, \mathbf{w})=\left(x_1, \ldots, x_t, y_{t+1}, \ldots, y_r, z, \mathbf{w}\right)$ centred at $q$ such that
$$
\begin{aligned}
\mathbf{u} & =\mathbf{x}^{\mathbf{A}_1}(\mathbf{y}-\boldsymbol{\xi})^{\boldsymbol{\Lambda}} \\
v & =\mathbf{x}^{\boldsymbol{\alpha}_1}(z-\zeta) \\
\mathbf{w} & =\mathbf{w},
\end{aligned}
$$

where $\zeta \neq 0, \xi_j \neq 0$ for all $j$ and the matrix $\boldsymbol{\Lambda}=\left(\lambda_{i, j}\right)$ of exponents has maximal rank, with $\lambda_{i, j} \in \mathbb{K}$. Moreover, if $\mathbf{u}^\gamma$ is not a first integral of $\omega \cdot \mathcal{O}_p$, then its total transform $\mathbf{u}^\gamma=\mathbf{x}^{\tilde{\gamma}} U$, where $U$ is a unit, satisfies one of the following:

- either the monomial $\mathbf{x}^{\widetilde{\gamma}}$ is not a first integral of $\widetilde{\omega}$;

- or, there exists a regular vector field $\partial_{y_i} \in \widetilde{\omega} \cdot \mathcal{O}_q$ such that $\partial_{y_i} U$ is a unit.
\end{lemma}

\begin{lemma}[Case 2, see {\cite[Lemma 10.3]{da2018local}}]\label{Lem:c2}
Assume that $\mathbf{A}_1$ does not have maximal rank. Then there exists a monomial coordinate system

$(\mathbf{x}, \mathbf{y}, \mathbf{w})=\left(x_1, \ldots, x_t, y_{t+1}, \ldots, y_{r+1}, \mathbf{w}\right)$ centred at $q$ such that
$$
\begin{aligned}
\mathbf{u} & =\mathbf{x}^{\mathbf{A}_1}(\mathbf{y}-\boldsymbol{\xi})^{\boldsymbol{\Lambda}} \\
v & =\mathbf{x}^{\boldsymbol{\alpha}_1} \\
\mathbf{w} & =\mathbf{w},
\end{aligned}
$$

where $\xi_j \neq 0$ for all $j$, the matrix $\mathbf{\Lambda}=\left(\lambda_{i, j}\right)$ is of maximal rank with entries in $\mathbb{K}$, and $\boldsymbol{\alpha}_1$ does not belong to the span of the rows of $\mathbf{A}_1$. Moreover, if $\mathbf{u}^\gamma$ is not a first integral of $\omega \cdot \mathcal{O}_p$, then its total transform $\mathbf{u}^\gamma=\mathbf{x}^{\widetilde{\gamma}} U$, where $U$ is a unit, satisfies one of the following:

- either the monomial $\mathbf{x}^{\widetilde{\gamma}}$ is not a first integral of $\widetilde{\omega}$;

- or there exists a regular vector field $\partial_{y_i} \in \widetilde{\omega} \cdot \mathcal{O}_q$ such that $\partial_{y_i} U$ is a unit.
\end{lemma}

In order to decrease $d_1$, we consider the ideal at $\pa$:
$$\mathcal{J}=\left(v^{d_1},\{\bu^{\bga_{2j}}v^j\}_{1\leq j\leq d_1-2},\{\bu^{\bga_{ij}}v^j\}_{3\leq i\leq N, 0\leq j\leq d_1-1},\bu^{\bbe}\right)$$

We now principalise $\mathcal{J}$ by a sequence of blow-ups $\tau$. Since $\mathcal{J}$ is generated by monomials of $\bu$ and $v$, the blow-up sequence $\tau$ can be chosen as combinatorial with respect to $\left(u_1\ldots u_rv_1=0\right)$. Let $\pb\in\tau^{-1}\left(\pa\right)$. We shall prove that the invariant at $\pb$ is strictly less than the invariant at $\pa$, which is $d_1$.

We distinguish two cases: Case 1 (Lemma \ref{Lem:c1}) and Case 2 (Lemma \ref{Lem:c2}), each of which is divided into three subcases, depending on which generator of $\mathcal{J}$ pulls back to generate $\tau^*\mathcal{J}$.

For Case 1, we have:
$$
\left\{
\begin{aligned}
\tau^* T_2&=U \mathbf{x}^{\mathbf{S}_{d_1}}(z-\zeta)^{d_1}+\sum_{j=1}^{d_1-2} \mathbf{x}^{\mathbf{S}_{ij}}(z-\zeta)^{j} c_{2j}(\mathbf{x}, \mathbf{y}, \mathbf{w})+\tau^*a_{20}(\mathbf{x}, \mathbf{y}, \mathbf{w}), \\
\tau^*T_{i_0} & = \tau^*\tilde{T}_{i_0}\mathbf{x}^{\mathbf{S}_{d_1}}(z-\zeta)^{d_1}+\sum_{j=1}^{d_1-1}\mathbf{x}^{\mathbf{S}_{i_0j}}(z-\zeta)^{j}c_{i_0j}+\mathbf{x}^{\mathbf{S}_{i_00}}c_{i_00}(\mathbf{x}, \mathbf{y}, \mathbf{w}),i\in \{3,\ldots,N\},
\end{aligned}
\right.
$$
where
\begin{enumerate}
\item the section $U$ is a unit which can be written as $\widetilde{U}(\mathbf{x}, \mathbf{y}, \mathbf{w})+\mathbf{x}^{\boldsymbol{\alpha}_1} \Omega(\mathbf{x}, \mathbf{y}, z, \mathbf{w})$, where $\widetilde{U}(\mathbf{x}, \mathbf{y}, \mathbf{w})$ is a unit and $\boldsymbol{\alpha}_1 \neq 0$ (because $\pb$ is in the pre-image of $\pa$);
\item for $i=2,\ldots,N,j=1,\ldots,d_1-1$, the section $c_{ij}$, pullback of $b_{ij}\left(\bu,\hat{\bv}\right)$, is either zero or a unit that does not depend on $z$;
\item the term $\mathbf{x}^{\mathbf{S}_{i_00}} c_{i_00}$ is the pullback of $a_{i_00}$; here either $c_{i_00}=0$ or it is equal to $w_m^\epsilon \widetilde{c}_{i_00}$, where $\widetilde{c}_{i_00}$ is a unit. Similarly, we write $\tau^*a_{i0}=\mathbf{x}^{\mathbf{S}_{i0}} c_{i0}$ where $c_{i0}$ is either zero or a unit that does not depend on $z$, for $i\neq i_0$.
\end{enumerate}

\begin{itemize}
\item[Case 1.1]{The pullback of $v^{d_1}$ generates $\tau^*\mathcal{J}$, i.e., $\mathbf{S}_{d_1}=\min \left\{\mathbf{S}_{d_1}, \mathbf{S}_{ij}, \mathbf{S}_{i_00}\right\}$.}
We have:
$$
\begin{aligned}
\tau^* T_2=\mathbf{x}^{\mathbf{S}_{d_1}} & {\left[\left(\widetilde{U} z+\widetilde{U} \zeta \nu+\mathbf{x}^{\boldsymbol{\alpha}_1} \Omega_2\right) z^{d_1-1}+\right.}\text { terms with power of } z<d_1-1],
\end{aligned}
$$
where $\boldsymbol{\alpha}_1$ is a non-zero matrix and $\Omega_2=[z+\zeta \nu] \Omega$. Since $\widetilde{U} z+\widetilde{U} \nu \zeta+\mathbf{x}^{\alpha_1} \Omega_2$ is a unit, we conclude that the invariant at $\pb\leq d_1-1<d_1$.

\item[Case 1.2]{There is a maximum $2\leq i_1\neq i_0\leq N,0\leq j_1\leq d_1-1$ or $i_1=i_0,1\leq j_1\leq d_1-1$ such that the pullback of $\mathbf{u}^{\bga_{i_1j_1}} v^{j_1}$ generates $\tau^*\mathcal{J}$, i.e., $\mathbf{S}_{i_1j_1}=\min \left\{\mathbf{S}_{d_1}, \mathbf{S}_{ij}, \mathbf{S}_{i_00}\right\}, \mathbf{S}_{d_1}>\mathbf{S}_{i_1j_1}$ and $\mathbf{S}_{i_1j}>\mathbf{S}_{i_1j_1}$ for $j>j_1$.}
We have:
$$
\tau^* T_{i_1}=\mathbf{x}^{\mathbf{S}_{i_1j_1}}\left[(z-\zeta)^{j_1} c_{i_1j_1}+\sum_{j=0}^{j_1-1} \mathbf{x}^{\mathbf{S}_{i_1j}-\mathbf{S}_{i_1j_1}}(z-\zeta)^j c_{i_1j}+\Psi(\mathbf{x}, \mathbf{y}, z, \mathbf{w})\right],
$$
where $\Psi(\mathbf{0}, \mathbf{y}, z, \mathbf{w}) \equiv 0$. Since $c_{i_1j_1}$ is a unit, we conclude that the invariant at $\pb\leq j_1<d_1$.

\item[Case 1.3]{The pullback of $\mathbf{u}^{\boldsymbol{\beta}}$ is the only generator of $\mathcal{J}^*$, that is, $\mathbf{S}_{i_00}=$ $\min \left\{\mathbf{S}_{d_1}, \mathbf{S}_{ij}, \mathbf{S}_{i_00}\right\}$ and $\mathbf{S}_{d_1}>\mathbf{S}_{i_00}, \mathbf{S}_{ij}>\mathbf{S}_{i_00}$.}
In this subcase, we consider two situations depending on the value of $\epsilon$:

Case 1.3a: $\epsilon=1$. Then
$$
\tau^* T_{i_0}=\mathbf{x}^{\mathbf{S}_{i_00}}\left[w_m \widetilde{c}_{i_00}+\Psi(\mathbf{x}, \mathbf{y}, z, \mathbf{w})\right],
$$
where $\widetilde{c}_{20}$ is a unit and $\Psi(\mathbf{0}, \mathbf{y}, z, \mathbf{w}) \equiv 0$, we conclude that the invariant at $\pb\leq 1<d_1$.

Case 1.3b: $\epsilon=0$. In this case, note that $\left(\boldsymbol{\beta}+\boldsymbol{\delta}\right)\wedge\bal\neq 0$. Thus, Lemma \ref{Lem:c1} guarantees that the total transform $\mathbf{u}^{\boldsymbol{\beta}+\boldsymbol{\delta}}=$ $\mathbf{x}^{\mathbf{S}_{i_00}+\widetilde{\boldsymbol{\delta}}} \widetilde{W}$, where $\widetilde{W}=(\mathbf{y}-\boldsymbol{\gamma})^{(\boldsymbol{\delta}+\boldsymbol{\beta}) \boldsymbol{\Lambda}}$ is a unit, satisfies one of the following:

- either $\left(\mathbf{S}_{i_00}+\widetilde{\boldsymbol{\delta}}\right)\wedge\tilde{\bal}\neq 0$, which implies that
$$
\tau^*\left[\mathbf{u}^\delta T_{i_0}\right]=\mathbf{x}^{\mathbf{S}_{i_00}+\widetilde{\boldsymbol{\delta}}} U
$$
for some unit $U$. Then we conclude that the invariant at $\pb=0<d_1$;

- or $\left(\mathbf{S}_{i_00}+\widetilde{\boldsymbol{\delta}}\right)\wedge\tilde{\bal}=0$, then we write
$$
\tau^*\left[\mathbf{u}^\delta T_{i_0}\right]=\mathbf{x}^{\mathbf{S}_{i_00}+\widetilde{\boldsymbol{\delta}}} \widetilde{W}(0)+\mathbf{x}^{\mathbf{S}_{i_00}+\widetilde{\boldsymbol{\delta}}}[\widetilde{W}-\widetilde{W}(0)+\Psi(\mathbf{x}, \mathbf{y}, z, \mathbf{w})],
$$
where $\Psi(\mathbf{0}, \mathbf{y}, z, \mathbf{w}) \equiv 0$ and the monomial $\mathbf{x}^{\mathbf{S}_{i_00}+\widetilde{\boldsymbol{\delta}}}W(0)\wedge\tilde{\bal}=0$. The linearly independent part of $\mathbf{x}^{\mathbf{S}_{i_00}+\widetilde{\boldsymbol{\delta}}}[\widetilde{W}-\widetilde{W}(0)+\Psi(\mathbf{x}, \mathbf{y}, z, \mathbf{w})]$ has order 1, so we conclude that the invariant at $\pb\leq 1<d_1$.
\end{itemize}

For Case 2, we have:
$$
\left\{
\begin{aligned}
\tau^* T_2&=U \mathbf{x}^{\mathbf{S}_{d_1}}+\sum_{j=1}^{d_1-2} \mathbf{x}^{\mathbf{S}_{2j}} c_{2j}(\mathbf{x}, \mathbf{y}, \mathbf{w})+\tau^*a_{20}(\mathbf{x}, \mathbf{y}, \mathbf{w}), \\
\tau^*T_{i_0} & = \tau^*\tilde{T}_{i_0}\mathbf{x}^{\mathbf{S}_{d_1}}+\sum_{j=1}^{d_1-1}\mathbf{x}^{\mathbf{S}_{i_0j}}c_{i_0j}+\mathbf{x}^{\mathbf{S}_{i_00}}c_{i_00}(\mathbf{x}, \mathbf{y}, \mathbf{w}),i\in \{3,\ldots,N\},
\end{aligned}
\right.
$$
where
\begin{enumerate}
\item the section $U$ is a unit;
\item for $i=2,\ldots,N,j=1,\ldots,d_1-1$, the section $c_{ij}$, pullback of $b_{ij}\left(\bu,\hat{\bv}\right)$, is either zero or a unit that does not depend on $z$;
\item the term $\mathbf{x}^{\mathbf{S}_{i_00}} c_{i_00}$ is the pullback of $a_{i_00}$; here either $c_{i_00}=0$ or it is equal to $w_m^\epsilon \widetilde{c}_{i_00}$, where $\widetilde{c}_{i_00}$ is a unit. Similarly, we write $\tau^*a_{i0}=\mathbf{x}^{\mathbf{S}_{i0}} c_{i0}$ where $c_{i0}$ is either zero or a unit that does not depend on $z$, for $i\neq i_0$;
\item note that:
$$
\mathbf{S}_{d_1}=d_1\boldsymbol{\alpha}_1 \quad \text { and } \quad \mathbf{S}_{ij}=j \boldsymbol{\alpha}_1+\bga_{ij} \mathbf{A}_1
$$

Hence each $\mathbf{S}_{d_1}$ and $\mathbf{S}_{ij}$ is a sum of an element of the span of the rows of $\mathbf{A}_1$ and a different multiple of the $\boldsymbol{\alpha}_1$. Since $\boldsymbol{\alpha}_1$ is linearly independent with the rows of $\mathbf{A}_1$, this means that the exponents $\mathbf{S}_{d_1}$ and $\mathbf{S}_{ij}$ are all distinct. Therefore, all of the multi-indices $\mathbf{S}_{ij}$ must be different.
\end{enumerate}

\begin{itemize}
\item[Case 2.1]\label{Case:2.1}{The pullback of $v^{d_1}$ generates $\tau^*\mathcal{J}$, i.e., $\mathbf{S}_{d_1}=\min \left\{\mathbf{S}_{d_1}, \mathbf{S}_{ij}, \mathbf{S}_{i_00}\right\}$.}
We have:
$$
\tau^*\left[\mathbf{u}^\delta T_2\right]=\mathbf{x}^{\mathbf{S}_{d_1}+\delta \mathbf{A}_1}[\widetilde{U}+\Psi(\mathbf{x}, \mathbf{y}, z, \mathbf{w})],
$$
where $\widetilde{U}=U(\mathbf{y}-\boldsymbol{\xi})^{\boldsymbol{\delta} \boldsymbol{\Lambda}}$ is a unit and $\Psi(\mathbf{0}, \mathbf{y}, z, \mathbf{w}) \equiv 0$. Since $\left(\mathbf{S}_{d_1}+\boldsymbol{\delta} \mathbf{A}_1\right)\wedge\tilde{\bal}\neq 0$ (which follows from the fact that $\boldsymbol{\alpha}_1$ does not belong to the span of the rows of $\mathbf{A}_1$), we conclude that the invariant at $\pb=0<d_1$.

\item[Case 2.2]{There is a maximum $2\leq i_1\neq i_0\leq N,0\leq j_1\leq d_1-1$ or $i_1=i_0,1\leq j_1\leq d_1-1$ such that the pullback of $\mathbf{u}^{\bga_{i_1j_1}} v^{j_1}$ is a generator of $\tau^*\mathcal{J}$, i.e., $\mathbf{S}_{i_1j_1}=\min \left\{\mathbf{S}_{d_1}, \mathbf{S}_{ij}, \mathbf{S}_{i_00}\right\}$ and $\mathbf{S}_{i_1j}>\mathbf{S}_{i_1j_1}$ for $j<j_1$.}
We have:
$$
\tau^*\left[\mathbf{u}^{\boldsymbol{\delta}} T_{i_1}\right]=\mathbf{x}^{\mathbf{S}_{i_1j_1}+\boldsymbol{\delta} \mathbf{A}_1}\left[\widetilde{c}_{i_1j_1}+\Psi(\mathbf{x}, \mathbf{y}, z, \mathbf{w})\right],
$$
where $\widetilde{c}_{i_1j_1}=c_{i_1j_1}(\mathbf{y}-\boldsymbol{\xi})^{\boldsymbol{\delta} \boldsymbol{\Lambda}}$ is a unit and $\Psi(\mathbf{0}, \mathbf{y}, z, \mathbf{w}) \equiv 0$. Since $\left(\mathbf{S}_{i_1j_1}+\widetilde{\boldsymbol{\delta}}\right)\wedge\tilde{\bal}$ (which follows from the fact that $\boldsymbol{\alpha}_1$ does not belong to the span of the rows of $\mathbf{A}_1$), we conclude that the invariant at $\pb=0<d_1$.

\item[Case 2.3]{The pullback of $\mathbf{u}^{\boldsymbol{\beta}}$ is the generator of $\mathcal{J}^*$ i.e., $\mathbf{S}_{i_00}=\min \left\{\mathbf{S}_{d_1}, \mathbf{S}_{ij}, \mathbf{S}_{i_00}\right\}$.}
In this subcase, we distinguish two possibilities according to the value of $\epsilon$:

Case 2.3a: $\epsilon=1$. Then
$$
\tau^* T_{i_0}=\mathbf{x}^{\mathbf{S}_{i_00}}\left[w_m \widetilde{c}_{20}+\Psi(\mathbf{x}, \mathbf{y}, z, \mathbf{w})\right]
$$

where $\widetilde{c}_{20}$ is a unit and $\Psi(\mathbf{0}, \mathbf{y}, z, \mathbf{w}) \equiv 0$, we conclude that the invariant at $\pb\leq 1<d_1$.

Case $2.3b, \epsilon=0$. In this case, note that $\left(\boldsymbol{\beta}+\boldsymbol{\delta}\right)\wedge\bal\neq 0$. Thus, Lemma \ref{Lem:c2} guarantees that the total transform $\mathbf{u}^{\boldsymbol{\beta}+\boldsymbol{\delta}}=$ $\mathbf{x}^{\mathbf{S}_{i_00}+\delta \mathbf{A}_1} \widetilde{W}$, where $\widetilde{W}=(\mathbf{y}-\gamma)^{(\boldsymbol{\delta}+\boldsymbol{\beta}) \boldsymbol{\Lambda}}$ is a unit, satisfies one of the following conditions:

- either $\left(\mathbf{S}_{i_00}+\boldsymbol{\delta} \mathbf{A}_1\right)\wedge\tilde{\bal}$, which implies that
$$
\tau^*\left[\mathbf{u}^{\boldsymbol{\delta}} T_{i_0}\right]=\mathbf{x}^{\mathbf{S}_{i_00}+\boldsymbol{\delta} \mathbf{A}_1} U
$$
for some unit $U$. We conclude that the invariant at $\pb=0<d_1$.

- or $\left(\mathbf{S}_{i_00}+\boldsymbol{\delta} \mathbf{A}_1\right)\wedge\tilde{\bal}=0$, then we write
$$
\tau^*\left[\mathbf{u}^\delta T_{i_0}\right]=\mathbf{x}^{\mathbf{S}_{i_00}+\boldsymbol{\delta} \mathbf{A}_1} \widetilde{W}(0)+\mathbf{x}^{\mathbf{S}_{i_00}+\boldsymbol{\delta} \mathbf{A}_1}[\widetilde{W}-\widetilde{W}(0)+\Psi(\mathbf{x}, \mathbf{y}, z, \mathbf{w})],
$$
where $\Psi(\mathbf{0}, \mathbf{y}, z, \mathbf{w}) \equiv 0$ and the monomial $\mathbf{x}^{\mathbf{S}_{i_00}+\boldsymbol{\delta} \mathbf{A}_1} W(0)\wedge\tilde{\bal}=0$. We conclude that the invariant at $\pb\leq 1<d_1$.
\end{itemize}
\end{proof}

\noindent
\textbf{Acknowledgements.} During the preparation of this work, the author was supported by the project ``Plan d’investissements France 2030'', IDEX UP ANR-18-IDEX-0001 at Université Paris Cité. The associated laboratory is Institut de Mathématiques de Jussieu-Paris Rive Gauche. I sincerely thank my supervisor André Belotto for his invaluable help and insightful comments, which have greatly contributed to the improvement and completion of this paper.

\bibliographystyle{amsalpha}
\bibliography{references}

------------------------------------------------------------------------------

\end{document}